\documentclass[12pt]{amsart}
\newcommand{\excise}[1]{}

\newtheorem{thm}{Theorem}[section]
\newtheorem{lemma}[thm]{Lemma}

\newtheorem{cor}[thm]{Corollary}
\newtheorem{prop}[thm]{Proposition}

\newtheorem{problem}[thm]{Problem}

\newtheorem{example}[thm]{Example}
\newtheorem{remark}[thm]{Remark}
\newtheorem{alg}[thm]{Algorithm}
\newtheorem{defn}[thm]{Definition}

\newenvironment{numbered}%
	{\begin{list}
		{\noindent\makebox[0mm][r]{\arabic{enumi}.}}
		{\leftmargin=5.5ex \usecounter{enumi}}
	}
	{\end{list}}

        {\begin{list}
                {\noindent\makebox[0mm][r]{(\roman{enumi})}}
                {\leftmargin=5.5ex \usecounter{enumi}}
        }
        {\end{list}}

\newenvironment{alglist}%
	{\begin{list}
		{}
		{\leftmargin=5.4em\labelwidth=5.0em\labelsep=.6em
		 \topsep=-1ex\itemsep=.1ex}\sf}  
	{\vspace{1ex}\end{list}}
\def\routine#1{\item[{\sc{#1}{\ }}]}
\def\procedure#1{{\sc{#1}}}
\newenvironment{routinelist}[1]%
	{\routine{#1}\begin{list}
		{}
		{\leftmargin=2.9em\labelwidth=2.4em\labelsep=.5em
		 \topsep=-1ex\itemsep=.1ex}}
	{\end{list}}
\newenvironment{algsublist}%
	{\begin{list}
		{}
		{\leftmargin=3.0em\labelwidth=4.8em\labelsep=.5em
		 \itemsep=.1ex\topsep=0ex}}
	{\end{list}}

\newenvironment{rcgraph}{\begin{trivlist}\item\centering\footnotesize$}
			{$\end{trivlist}}

\def\s{\scriptstyle}
\def\sst{\scriptscriptstyle}
\def\dis{\displaystyle}
\def\fillrightmap{\mathord- \mkern-6mu
	\cleaders\hbox{$\mkern-2mu \mathord- \mkern-2mu$}\hfill
	\mkern-6mu \mathord\rightarrow}
\newcommand{\injmatrix}[3]{{
\begin{array}{@{}r@{\;}r@{}c@{}l@{}}
  \begin{array}{@{}c@{}}		
	\begin{array}{@{}c@{\ }c@{}}
	\\\\
	#1
	\end{array}\!
  \end{array}
&
  \begin{array}{@{}c@{}}		
	\begin{array}{@{}l@{}}\\\\[.5ex]			
	\end{array}						
	\\							
	\left(\begin{array}{@{}l@{}}				
	#3							
	\end{array}\!						
	\right.							
  \end{array}							
&
  #2					
&
  \begin{array}{@{}c@{}}		
	\begin{array}{@{}l@{}}\\\\				
	\end{array}						
	\\							
	\left.\!\begin{array}{@{}l@{}}				
	#3							
	\end{array}						
	\right)							
  \end{array}							
\end{array}
}}

\def\bem#1{\textbf{#1}}

\def\<{\langle}
\def\>{\rangle}
\def\0{{\mathbf 0}}
\def\1{{\mathbf 1}}
\def\GF{\Gamma_{\hspace{-.4ex}F\hspace{.1ex}}}
\def\GI{\Gamma_{\hspace{-.4ex}I\hspace{.1ex}}}
\def\Gm{\Gamma_{\hspace{-.4ex}\mm\hspace{.1ex}}}
\def\AA{{\mathcal A}}

\def\FF{{\mathcal F}}

\def\MM{{\mathcal M}}
\def\NN{{\mathbb N}}
\def\RR{{\mathbb R}}
\def\SS{{\mathcal S}}
\def\WW{{\hspace{.1ex}\overline{\hspace{-.1ex}W\hspace{-.1ex}}\hspace{.1ex}}{}}
\def\ZZ{{\mathbb Z}}
\def\aa{\alpha}
\def\bb{\beta}
\DeclareMathSymbol{\Bbbk}{\mathord}{AMSb}{"7C}
\def\kk{\Bbbk}

\def\mm{{\mathfrak m}}
\def\xx{{\mathbf x}}

\def\th{{\rm th}}
\def\im{{\rm im}}

\def\sat{{\rm sat}}
\def\too{\longrightarrow}

\def\into{\hookrightarrow}
\def\onto{\to}
\def\spot{{\hbox{\raisebox{1.7pt}{\large\bf .}}\hspace{-.5pt}}}
\def\minus{\setminus}

\def\nothing{\emptyset}
\def\bigcupdot{\makebox[0pt][l]{$\hspace{1.05ex}\cdot$}\textstyle\bigcup}

\def\ol#1{{\overline {#1}}}

\def\mod#1{\ ({\rm mod}\:#1)}

\begin{document}

\raggedbottom

\begin{excise}{
\title[Algorithms for graded injective resolutions and local cohomology]
{Algorithms for graded injective resolutions\\ and local cohomology over
semigroup rings}
\author[David Helm and Ezra Miller]
	{David Helm and Ezra Miller\\
	Mathematics Department, Harvard University, Cambridge, MA, USA
	\nextaddress
	Mathematical Sciences Research Institute, Berkeley, CA, USA
	\nextaddress (current address:
	School of Mathematics, University of Minnesota, Minneapolis, MN, USA)
	}
\maketitle

\begin{abstract}
%
Let $Q$ be an affine semigroup generating~${\s \ZZ}^d$, and fix a
finitely generated ${\s \ZZ}^d$-graded module~$M$ over the semigroup
algebra~${\kk}[Q]$ for a field~${\kk}$.  We provide an algorithm to
compute a minimal ${\s \ZZ}^d$-graded injective resolution of~$M$ up to
any desired cohomological degree.  As an application, we derive an
algorithm computing the local cohomology modules~$H^i_I(M)$ supported
on any monomial (that is, ${\s \ZZ}^d$-graded) ideal~$I$.
Since
these local cohomology modules are neither finitely generated nor
finitely cogenerated, part of this task is defining a finite data
structure to encode~them.
\end{abstract}
}\end{excise}

\title[Algorithms for graded injective resolutions and local cohomology]
{Algorithms for graded injective resolutions\\ and local cohomology over
semigroup rings}
\author{David Helm}
\address{Mathematics Department\\UC Berkeley\\Berkeley, CA}
\email{dhelm@math.harvard.edu}
\author{Ezra Miller}
\address{Mathematical Sciences Research Institute\\Berkeley, CA}
\email{ezra@math.umn.edu}
\raisebox{0ex}{\date{15 September 2003}}
\begin{abstract}
%
Let $Q$ be an affine semigroup generating~${\ZZ}^d$, and fix a
finitely generated ${\ZZ}^d$-graded module~$M$ over the semigroup
algebra~${\kk}[Q]$ for a field~${\kk}$.  We provide an algorithm to
compute a minimal ${\ZZ}^d$-graded injective resolution of~$M$ up to
any desired cohomological degree.  As an application, we derive an
algorithm computing the local cohomology modules~$H^i_I(M)$ supported
on any monomial (that is, ${\ZZ}^d$-graded) ideal~$I$.
Since
these local cohomology modules are neither finitely generated nor
finitely cogenerated, part of this task is defining a finite data
structure to encode~them.
\end{abstract}
\maketitle

{}

\setcounter{tocdepth}{1}
\tableofcontents

\section{Introduction}

\label{sec:intro}

Injective resolutions are fundamental homological objects in
commutative algebra.  For general noetherian rings with arbitrary
gradings, however, injective modules are so big, and injective
resolutions so intractable, that effective computations are never made
using them.  But when the ring in question is an affine semigroup ring
of dimension~$d$, the natural grading by~$\ZZ^d$ is substantially
better behaved: $\ZZ^d$-graded injective modules can be expressed
polyhedrally and are therefore quite explicit.  In this paper we
provide algorithms to compute $\ZZ^d$-graded injective resolutions
over affine semigroups rings.  Part of this task is finding a finite
data structure to express the output.

As an application, we provide an algorithm to compute the local
cohomology, supported on an arbitrary monomial ideal, of a finitely
generated $\ZZ^d$-graded module over a normal affine semigroup ring.
As far as we are aware, this is the first algorithm to compute
local cohomology for any general class of modules over any class of
nonregular rings.
Our motivation was to make a systematic study of conditions on the
support ideal and the ambient ring that cause local cohomology to have
infinite Bass numbers.  That such infinite behavior occurs only over
nonregular rings necessitated our working over affine semigroup rings,
which seem to be the simplest available singular rings.

To make our context precise, let $Q \subset \ZZ^d$ be an affine
semigroup, that is, a finitely generated submonoid of~$\ZZ^d$.  We
assume that $Q$ is sharp, meaning that $Q$~has no units, and that $Q$
generates~$\ZZ^d$ as a group.  Consider the semigroup algebra $\kk[Q]
= \bigoplus_{a \in Q} \kk\cdot\{\xx^a\}$ over a field~$\kk$.  The
modules that concern us comprise the category~$\MM$ of $\ZZ^d$-graded
modules $H = \bigoplus_{\alpha\in\raisebox{-.15ex}{${\sst \ZZ}$}^d}
H_\alpha$ for which there exists a bound independent of~$\alpha$ on
the dimensions of the graded pieces~$H_\alpha$ as vector spaces
over~$\kk$.  The injective objects in~$\MM$ are described in
Section~\ref{sec:effective}, and every finitely generated
$\ZZ^d$-graded module lies in~$\MM$.  Our main theorem concerning
injectives is the following.

\begin{thm} \label{t:injres'}
Fix a finitely generated $\ZZ^d$-graded module $M$ over an affine
semigroup ring~$\kk[Q]$ and an integer~\mbox{$n \geq 0$}.  The first
$n$ stages in a minimal $\ZZ^d$-graded injective resolution of~$M$ can
be expressed in a finite, algorithmically computable data structure.
\end{thm}

A more precise version, along with a pointer to the algorithms that do
the job, is stated in Theorem~\ref{t:injres}.  The data structure
consists of a list of {\em monomial matrices}, as we define in
Section~\ref{sec:effective}, generalizing those for $Q = \NN^d$
in~\cite{Mil2}.  The idea of the algorithm in Theorem~\ref{t:injres'}
is to do all computations using {\em irreducible resolutions}\/
\cite{MilIrr} as faithful approximations to injective resolutions.
Background on {\em irreducible hulls}\/ is presented in
Section~\ref{sec:effective}; the algorithms for working with them
constitute Section~\ref{sec:irred}.  The derivation of an algorithm
for injective resolutions is then completed in
Section~\ref{sec:injres}.

Even more seriously than is the case with injective resolutions, a
substantial part of building an algorithm to compute local cohomology
is finding a finite data structure to express the output.  Indeed,
unlike injectives in our category~$\MM$, and in stark contrast with
the regular case (even without a grading \cite{HuSh,Lyu1}), the local
cohomology $H^i_I(M)$ often has neither a finite generating set {\em
nor a finite cogenerating set}\/ \cite{HarCofinite,HelM}.  This
remains true even when $M$ is finitely generated and $I \subseteq
\kk[Q]$ is a $\ZZ^d$-graded ideal---that is, generated by monomials.
\begin{excise}{%
  This raises an important combinatorial question, which remains open
  even for saturated semigroups~$Q$.
  
  \begin{problem} \label{prob:socle}
  Characterize the monomial ideals $I \subset \kk[Q]$ such that the
  local cohomology $H^i_I(\kk[Q])$ of
  the ring\/~$\kk[Q]$ supported on~$I$ has infinite-dimensional
  socle.
  \end{problem}
  
  In other words, when does the local cohomology contain an
  infinite-dimensional vector subspace annihilated by the maximal
  $\ZZ^d$-graded ideal of~$\kk[Q]$, which is generated by all nonunit
  monomials?  One could ask about infinite zeroth Bass numbers at
  other $\ZZ^d$-graded primes, but even understanding the socle would
  be a big step.  The issue is combinatorial because we may assume $I$
  is radical, and such monomial ideals correspond to
  polyhedral subcomplexes of the real cone generated by~$Q$.
  
  Problem~\ref{prob:socle} remains open in part because there are no
  known combinatorial descriptions even of the Hilbert series for
  local cohomology of~$\kk[Q]$ with arbitrary support, let alone its
  module structure.  There are ``minimal'' complexes of localizations
  that whose homology is the local cohomology $H^i_I(\kk[Q])$
  \cite{YanMonSup}, but even these {\em canonical \v Cech complexes}\/
  require the equivalent of a computation of a minimal injective
  resolution of a monomial quotient of a normal semigroup ring.
}\end{excise}%
Our solution is to
decompose~$\ZZ^d$ into tractable regions on which the local cohomology
is constant.

\begin{defn}\rm \label{d:sector}
Suppose $H$ is a $\ZZ^d$-graded module over an affine semigroup
ring~$\kk[Q]$.  A \bem{sector partition} of~$H$ is
\begin{numbered}
\item
a finite partition $\ZZ^d = \bigcupdot_{S \in \SS} S$ of the lattice
$\ZZ^d$ into \bem{sectors}, each of which is required to consist of
the lattice points in a finite disjoint union of rational polyhedra
defined as intersections of half-spaces for hyperplanes
parallel~to~facets~of~$Q$;

\item
a finite-dimensional vector space $H_S$ for each sector $S \in \SS$,
along with isomorphisms $H_\aa \to H_S$ for all $\ZZ^d$-graded degrees
$\aa \in S$; and

\item
vector space homomorphisms $H_S \stackrel{\xx^{T-S}}\too H_T$ whenever
there exist $\aa \in S$ and $\bb \in T$ satisfying $\bb-\aa \in Q$,
such that for all choices of $\alpha$ and~$\beta$, the
diagram~commutes:
$$
\begin{array}{c@{\ }c@{\ }c}
           & \scriptstyle\xx^{\bb-\aa} &            \\[-1.5ex]
    H_\aa  &      \fillrightmap        &    H_\bb   \\
\downarrow &        \ \qquad\          & \downarrow \\[-1ex]
           &   \scriptstyle\xx^{T-S}   &            \\[-1.5ex]
     H_S   &      \fillrightmap        &     H_T
\end{array}
$$
\end{numbered}
Write $\SS \vdash H$ to indicate the above sector partition.  (The
commutativity of the above diagram implies immediately that
$\xx^{S-S}$ is the identity, and that $\xx^{R-T}\xx^{T-S} =
\xx^{R-S}$.)
\end{defn}

%
%
The finite data structure of a sector partition $\SS \vdash H$,
including the spaces~$H_S$ and the maps~$\xx^{T-S}$, clearly suffice
to reconstruct~$H$ up to isomorphism.
%
%
The second half of this paper is devoted to
computing sector partitions for~$H$ when $H = H^i_I(M)$ is a local
cohomology module.

\begin{thm} \label{t:sector}
For any finitely generated $\ZZ^d$-graded module~$M$ over a normal
semigroup ring~$\kk[Q]$ and any monomial ideal~$I$, each local
cohomology module~$H^i_I(M)$ has an algorithmically computable sector
partition $\SS \vdash H^i_I(M)$.
\end{thm}

Section~\ref{sec:decomp} demonstrates how sector partitions arise for
the cohomology of any complex of injectives over a normal semigroup
ring.  Algorithms for producing these sector partitions, particlularly
the expressions of sectors as unions of polyhedral sets of lattice
points, occupy Section~\ref{sec:secposet}.  The proof of
Theorem~\ref{t:sector}, by expressing local cohomology as the
cohomology of a complex of injectives (algorithmically computed by
Theorem~\ref{t:injres'}) in the usual way, occurs in
Section~\ref{sec:injalg}.  That section also treats complexity issues.
The main thrust is that for fixed dimension~$d$, the running times of
our algorithms are all polynomial in the Bass numbers of the finitely
generated input module~$M$ and the number of facets of~$Q$, times the
usual factor arising from the complexity of Gr\"obner basis
computation, where it occurs.  If~$d$ is allowed to vary, then the
numbers of polyhedra comprising sectors increase exponentially
with~$d$.


Theorem~\ref{t:sector} allows the computation of many features of
local cohomology modules.  For example, Hilbert series simply record
the vector space dimensions in each of the finitely many sectors $S
\in \SS$.  Our algorithms can actually calculate these dimensions
without computing the maps in part~3 of Definition~\ref{d:sector},
making it easier to determine when (for example) $H^i_I(M)$ is
nonzero.  Future algorithmic methods (currently open problems) include
the calculation of associated primes and locations of socle degrees
(even if there are infinitely many)
using a sector partition as input.  In particular, because of the
finiteness of the number of polyhedra partitioning sectors, we believe
that the socle degrees should lie along polyhedrally describable
subsets of~$\ZZ^d$.

\subsection*{Historical context}

There have been a number of recent algorithmic computations in local
cohomology, such as those by Walther \cite{Wal} (based on abstract
methods of Lyubeznik \cite{Lyu1}),
Eisenbud--Musta\c{t}\v{a}--Stillman \cite{EMS}, Miller \cite{Mil2},
Musta\c{t}\v{a} \cite{Mus1}, and Yanagawa \cite{YanMonSup}.  These and
related papers fall naturally into a number of categories.  For
instance, the last three deal with $\ZZ^d$-graded modules over
polynomial rings in $d$ variables; in particular, they compute local
cohomology with support on monomial ideals.  In contrast, the paper
\cite{EMS} works with coarser gradings---but still
with monomial support, while \cite{Wal} requires no grading at all.
As the gradings used become coarser, the papers increasingly depend on
Gr\"obner bases: the monomial ideal papers require very little (if
any) Gr\"obner basis computation; the coarser gradings depend heavily
on commutative Gr\"obner bases; and the nongraded methods rely on
noncommutative Gr\"obner bases over the Weyl algebra.

Regardless of the methods, all of the above papers share one
fundamental aspect: the base ring is regular (usually a polynomial
ring, in the algorithmic setting).  The reason for restricting to
these rings
is that local cohomology over them behaves in many respects like a
finitely generated module, even though it usually fails to be finitely
generated.  For example, Lyubeznik \cite{Lyu1} and Walther \cite{Wal}
take advantage of the fact that local cohomology modules over regular
rings are finitely generated (indeed, holonomic) over the
corresponding algebra of differential operators, and that the algebra
of differential operators of a regular ring is easily presented, at
least in characteristic zero.

Generally speaking, our methods lie somewhere between the monomial and
coarsely graded methods described above, relying on a mix of Gr\"obner
bases and integer programming.  The principle underlying our
computation of injective resolutions is that one should attempt to
recover entire $\ZZ^d$-graded modules from their $Q$-graded parts.
This idea originated for polynomial rings in
\cite{Mil9812095,Mus1,Mil2}, was transfered in a restricted form to
semigroup rings in \cite{YanPoset}, and developed generally for
semigroup-graded noetherian rings in \cite{HelM}.  In the present
context, the recovery of a module from its $Q$-graded part suggested
that we compute injective resolutions via the irreducible resolutions
of \cite{MilIrr}.

Origins of the notion of sector partition can be seen in the Hilbert
series formula for the local cohomology of canonical modules of normal
semigroup rings \cite{TerLocalCoh,YanMonSup}, where the cellular
homology was constant on large polyhedral regions of~$\ZZ^d$.  The
accompanying notion of {\em straight module}\/ \cite{YanPoset,HelM}
abstracted this constancy; in fact, our Theorem~\ref{t:decomp} is
really a theorem about straight modules as in
\cite[Definition~5.1]{HelM}.  In any case, once the injective
resolution has been computed using irreducible resolutions, the sector
partition for local cohomology requires the entire $\ZZ^d$-graded
structure of the injective resolution, and not just its $Q$-graded
part.

\begin{excise}{%
  \begin{remark}\rm
  Only in the regular case are the injective resolutions
  computed here actually equivalent by Alexander duality to free
  resolutions \cite{Mil2}.  This follows essentially because in
  nonregular semigroup rings there exist pairs of principal ideals
  whose intersections fail to be principal \cite{HelM}.
  \end{remark}

  Commutative Gr\"obner bases pervade our algorithms, even if only
  implicitly (other methods for making the standard computations could
  be substituted).  Their role is heavier than in the case $Q =
  \NN^d$, where injective resolutions can be computed as the Alexander
  duals to free resolutions \cite[Section~4.2]{Mil2}.  For nonregular
  semigroup rings, free and injective resolutions reflect distinctly
  different homological and polyhedral information.  which comes down
  to the fact that the set of lattice points in a rational translate
  of the real cone $-\RR_+Q$ need not have a unique maximal element.
}\end{excise}%

\subsection*{Conventions and notation}
%
In addition to the notation introduced thus far, we close this
Introduction with a note on conventions.  The semigroup~$Q$ is
required to be saturated in Sections~\ref{sec:decomp}--%
\ref{sec:injalg} because we do not know how to compute sector
partitions in the unsaturated context (Remark~\ref{rk:holes}).  Other
than the temporary saturation requirement in
Section~\ref{sub:satgenrel}, the semigroup can be unsaturated in
Sections~\ref{sec:effective}--\ref{sec:injres}.
(Reminders of these conventions appear in each section).

The symbol $\xx^\alpha \in \kk[\ZZ^d]$ denotes a Laurent monomial in
the localization $\kk[\ZZ^d]$ of the semigroup ring~$\kk[Q]$.  The
$\kk$-vector space spanned by $\{\xx^\alpha \mid \alpha \in T\}$ for a
subset $T \subseteq \ZZ^d$ will be denoted by $\kk\{T\}$.  The
$\kk$-subalgebra of $\kk[\ZZ^d]$ will be denoted by~$\kk[T]$.

The \bem{faces} of $Q$ are those subsets minimizing linear functionals
on~$Q$.  The \bem{edges} and \bem{facets} are the faces of
dimension~$1$ and codimension~$1$.  To every face~$F$ corresponds a
prime ideal~$P_F$ and a quotient affine semigroup ring $\kk[F] =
\kk\{F\}$.

All modules in this paper are $\ZZ^d$-graded unless otherwise stated.
In particular, injective modules (defined in Section~\ref{sec:decomp})
are $\ZZ^d$-graded injective, which means that they are usually not
injective in the category of all $\kk[Q]$-modules.  Two subsets $S,T
\subseteq \ZZ^d$ have the difference set $T-S = \{\beta-\alpha \mid
\alpha \in S \hbox{\rm\ and } \beta \in T\} \subseteq \ZZ^d$.  This
allows us to write the localization of~$M$ along a face~$F$ as the
module $M[\ZZ F] := M \otimes_{\kk[Q]} \kk[Q-F]$.
Homomorphisms $N \to N'$ of modules are assumed to have $\ZZ^d$-graded
degree~$\0$, so that $N_\alpha \to N'{}_{\hspace{-.7ex}\alpha}$ for
all $\alpha \in \ZZ^d$.

We assume in this paper that standard algorithmic calculations with
finitely generated modules over~$\kk[Q]$ are available.  In
particular, we assume that the homology of any three-term (nonexact)
sequence of finitely generated modules can be calculated, as can the
submodule annihilated by a prime ideal of~$\kk[Q]$.  The
$\ZZ^d$-grading only makes these computations easier, and the results
of all such algorithms are still $\ZZ^d$-graded.

\section{Effective irreducible hulls}

\label{sec:effective}

In this section the affine semigroup~$Q$ need not be saturated.  In
the $\ZZ^d$-graded category~$\MM$ from the Introduction, the injective
modules have simple~descriptions.

\begin{defn}\rm \label{d:inj}
Let $T \subset \ZZ^d$ be closed under addition of elements of~$-Q$, by
which we mean $T - Q \subset T$.  Then $\kk\{T\}$ can be given the
structure of a $\kk[Q]$-module by setting
\begin{eqnarray*} \label{eq:T}
  \xx^a \xx^\beta \ \:=\ \:
	\left\{\begin{array}{@{\ }ll}
		\xx^{a + \beta}& \hbox{if } a + \beta \in T\\
		0		& \hbox{otherwise}.
	\end{array}\right.
\end{eqnarray*}
An \bem{indecomposable injective} is any module of the form
$\kk\{\alpha + F-Q\}$, for some face $F$ and $\alpha \in \ZZ^d$.  
\end{defn}

All such objects are injective in~$\MM$, and every injective object of
$\MM$ is isomorphic to a finite direct sum of indecomposable
injectives \cite[Chapter~11]{cca}.  We shall work exclusively with
objects in~$\MM$.  Thus the term ``injective module''
in the rest of this paper will refer to modules of the above type.

Injectives are infinitely generated.  For computations, we therefore
work with certain finitely generated approximations.  A module $N$ is
called \bem{$Q$-graded} if $N$ equals its \bem{$Q$-graded part} $N_Q
:= \bigoplus_{a \in Q} N_a$.  A submodule $N$ of a module $N'$ is an
\bem{essential submodule} if $N$ intersects every nonzero submodule of
$N'$ nontrivially; the inclusion $N \into N'$ is also called an
\bem{essential extension}.  In particular, $N$ must be nonzero.

\begin{defn}\rm \label{d:irred}
An \bem{irreducible sum} is a module that can be expressed as the
$Q$-graded part $J_Q$ of some injective module~$J$.  An
\bem{irreducible hull} of a $Q$-graded module $N$ is an irreducible
sum $\WW$ along with an essential extension $N \into \WW$.
\end{defn}

The existence of unique minimal injective resolutions
\cite[Corollary~11.35]{cca} includes the fact that every finitely
generated module has an injective hull (that is, an inclusion into an
injective that is an essential extension) that is unique up to
isomorphism.  Taking $Q$-graded parts yields immediately the following
lemma.

\begin{lemma} \label{l:irredHull}
Every $Q$-graded module has an irreducible hull.  It is unique up to
isomorphism, and isomorphic to the $Q$-graded part of an injective
hull of~$M$.
\end{lemma}

We call the modules of Definition~\ref{d:irred} irreducible sums
because of the next lemma, which is \cite[Lemma~2.2]{MilIrr}.  An
ideal $W$ is called \bem{irreducible} if
$W$ can not be expressed as an intersection of two ideals properly
containing~it.

\begin{lemma} \label{l:irrideal}
A monomial ideal $W$ is irreducible
if and only if the $Q$-graded part of some indecomposable injective
module~$J$ satisfies $J_Q = \WW$.
\end{lemma}

Modules $M$ are usually stored as data structures keeping track of
their generators and relations---that is, as quotients of free
modules.  In the context of injective resolutions and local
cohomology, storing $M$ as a submodule of an irreducible sum is also
useful.  Our next definition specifies a data structure that precisely
describes an irreducible sum~$\WW$.

\begin{defn}\rm \label{d:WW}
\bem{Effective data} for an irreducible sum $\WW = \bigoplus_{j=1}^r
\kk\{\alpha_j + F_j-Q\}_Q$ consist of:
\begin{numbered}
\item
an ordered $r$-tuple $F_1,\ldots,F_r$ of faces of $Q$; and

\item
an ordered $r$-tuple $\alpha_1,\ldots,\alpha_r$, where $\alpha_j \in
\ZZ^d/\ZZ F_j$ satisfies $Q \cap (\alpha_j + \ZZ F_j) \neq \nothing$.
\end{numbered}
An \bem{effective vector} of degree $a \in Q$ is an $r$-tuple
$(\lambda_1,\ldots,\lambda_r) \in \kk^r$ such that $\lambda_j = 0$
whenever $a \not\in \alpha_j + F_j-Q$.  Concatenation of the
respective face and degree data from two effective data yields their
\bem{direct sum}.
\end{defn}

Note that the faces $F_j$ need not be distinct.  The condition $\alpha
\in \ZZ^d/\ZZ F$ takes care of the fact that two degrees $\alpha$ and
$\alpha'$ off by an element of $\ZZ F$ give the same module
\hbox{$\kk\{\alpha+F-Q\} = \kk\{\alpha'+F-Q\}$}.  Usually the
$\alpha$'s are recorded as elements of~$\ZZ^d$, since the quotient mod
$\ZZ F$ can be deduced from the face data.  The condition \hbox{$Q
\cap (\alpha + \ZZ F) \neq \nothing$} ensures that $\kk\{\alpha +
F-Q\}$ has nonzero $Q$-graded part.  The condition on the $\lambda$'s
simply requires each nonzero component to lie in a nonzero degree of
the corresponding irreducible summand.

\begin{defn}\rm \label{d:MinWW}
An \bem{effective irreducible hull} of a $Q$-graded module~$M$ consists
of effective data for $\WW$ plus a list of finitely many effective
vectors in\/ $\WW$ generating a submodule isomorphic to~$M$.
\end{defn}

An irreducible hull $M \into \WW$ is not quite dual to an expression $\FF
\onto M$ as a quotient of a free module.  The generators of $\FF$ have as
their dual notion the face data $F_1,\ldots,F_r$, which as abstract
objects associated to $M$ are known as \bem{cogenerators}.  Just as the
degrees of the generators of $\FF$ need to specified, so must the degree
data for the cogenerators.  However, the notion of effective vector for
$M$ as a submodule of~$\WW$ is dual not to the notion of relation for
$M$ inside~$\FF$, but rather to the notion of cogenerator for~$M$.
Relations for $M$ are, in
actuality,
dual to the notion of \emph{co}generators for the \emph{co}kernel of
$M \into \WW$, which correspond to indecomposable summands in
cohomological degree~$1$ of the minimal injective resolution of~$M$;
we dub these the \bem{correlations} of~$M$.  Thus a presentation of
$M$ by generators and relations is dual to a presentation of $M$ by
cogenerators and correlations, whereas an irreducible hull presents
$M$ by generators and cogenerators.

\section{Computing with irreducible hulls}

\label{sec:irred}

Given a $Q$-graded module~$M$ in the usual way, via generators and
relations, this section computes an irreducible hull $M \into \WW$ as
well as the cokernel of this inclusion.

Calculating an effective irreducible hull of~$M$ is, by definition,
equivalent to calculating an irreducible decomposition of~$M$.
Thinking of the case $M = \kk[Q]/I$ for a monomial ideal~$I$, this
procedure is polyhedral in nature: it writes the set of monomials
outside of~$I$ as a union of convex polyhedral regions whose facets
are parallel to those of~$Q$.  The algorithm for computing an
effective irreducible hull $M \into \WW$, culminating in
Proposition~\ref{p:irred}, does not require $Q$ to be saturated.

Computing the cokernel, however, is strictly easier for saturated
semigroups.  The main point is the computation of generators for
irreducible ideals.  For saturated semigroups this is
Proposition~\ref{p:saturated}.  The harder unsaturated case, in
Proposition~\ref{p:W}, relies on the computation of irreducible ideals
over its saturation.  To highlight the simplification in the saturated
case, we state the main result of
Sections~\ref{sub:satgenrel}--\ref{sub:unsatgenrel}~here.

\begin{prop} \label{p:relations}
Generators and relations for~$M$ and\/ $\WW/M$ are algorithmically
computable from an effective irreducible hull $M \into \WW$ over any
affine semigroup ring\/~$\kk[Q]$.
\end{prop}
\begin{proof}
Generators for $\WW$ are already given, and relations for $\WW$
constitute a direct sum of irreducible ideals calculated as in
Proposition~\ref{p:saturated} for saturated semigroups, and
Proposition~\ref{p:W} in general.  Since $M$ is specified by its
generators as a submodule of~$\WW$, the current proposition reduces to
calculating submodules and quotients of modules presented by
generators and relations.
\end{proof}

\subsection{Effective irreducible hulls from generators and relations}


\label{sub:effIrr}

This subsection does not require the affine semigroup~$Q$ to be
saturated.  The next two results make Algorithm~\ref{a:irred} possible
to state and easier to read.  The notation $\<y_1,\ldots,y_j\>$ means
`the $\kk[Q]$-submodule generated by the elements $y_1,\ldots,y_j$ in
their ambient module', and $(0:_M P_F)$ is the submodule of $M$
annihilated by~$P_F$.

\begin{lemma} \label{l:B}
Suppose $F$ has minimal dimension among faces of $Q$ such that $P_F$
is associated to~$M$.  Then the natural map $(0:_M P_F)$ to its
localization $(0:_M P_F)[\ZZ F]$ along $F$ is an inclusion.
Furthermore, we can find algorithmically a set $B \subset (0:_M P_F)$
of homogeneous elements that consitute a $\kk[\ZZ F]$-basis for $(0:_M
P_F)[\ZZ F]$.
\end{lemma}
\begin{proof}
The $\kk[Q]$-module $(0:_M P_F)$ is naturally a torsion-free
$\kk[F]$-module, by minimality of $\dim F$.  Therefore $(0:_M P_F)$
includes into its localization along~$F$, which must be a free
$\kk[\ZZ F]$-module.  Now use the following algorithm.
\end{proof}

\begin{alg}[for Lemma~\ref{l:B}]\rm \label{a:B}
Choose any element of $(0:_M P_F)$ as the first basis vector~$y_1 \in
B$.  Having chosen~$y_j$, let $y_{j+1}$ be any element of $(0:_M P_F)$
whose image in $(0:_M P_F)/\<y_1,\ldots,y_j\>$ generates a submodule
of Krull dimension $\dim F$ (equivalently, the image of $y_{j+1}$ has
annihilator~$P_F$).  The algorithm terminates when the Krull dimension
of the quotient $(0:_M P_F)/\<y_1,\ldots,y_j\>$ is strictly less than
$\dim F$.
\end{alg}

\begin{lemma} \label{l:B'}
In the situation of Lemma~\ref{l:B}, the scalar factor on the
(monomial) coefficient of \hbox{$y \in B$} in the unique $\kk[\ZZ
F]$-linear combination of elements in~$B$ equaling any fixed element
\hbox{$z \in (0:_M P_F)$} can be computed algorithmically.
\end{lemma}

We present the proof as an algorithm.

\begin{alg}[for Lemma~\ref{l:B'}]\rm \label{a:B'}
Let $B(z) = \{y \in B \mid \deg(y) \equiv \deg(z) \mod{\ZZ F}\}$.  The
coefficient of $y$ in~$z$ is zero if $y \not\in B(z)$.  Otherwise,
find elements $a$ and $\{a_y \mid y \in B(z)\}$ in the face $F$ such
that $a + \deg(z) = a_y + \deg(y)$ for all $y \in B(z)$.  By
construction, $\{\xx^{a_y} \cdot y \mid y \in B(z)\}$ is a $\kk$-basis
for the degree $a + \deg(z)$ piece of $(0:_M P_F)$, and standard
methods allow us to calculate the syzygy with $\xx^a \cdot z$.
\end{alg}

Write $\GF N := \Gamma_{\hspace{-.4ex}P_F} N = (0:_N P_F^\infty)$ for
the set of elements in $N$ annihilated by all high powers of~$P_F$.

\begin{alg}\rm \label{a:irred}
\end{alg}
\begin{alglist}
\routine{input} $Q$-graded module $M$ given by a generating set
	$G \subset M$ and relations

\routine{output} effective irreducible hull $\WW$ of $M$ with
	effective vector set $\Lambda$ \hbox{indexed by $G$}

\begin{routinelist}{initialize}
\item[$N :=$] $M$

\item[$\WW :=$] $(\{\},\{\})$, the empty effective datum for the
	irreducible sum $0$

\item[$\lambda_g :=$] $()$ for all $g \in G$; here $()$ is the
	effective vector of length zero in $\WW$

\item[$i :=$] $1$
\end{routinelist}

\begin{routinelist}{define}
\item[] \hspace{-2.9em}$(F_1,\ldots,F_s) :=$ an ordering of the faces
	of $Q$ with $\dim(F_i) \leq \dim(F_{i+1})$
\end{routinelist}

\routine{while} $i \leq s$ \procedure{do}

\begin{algsublist}

	\begin{routinelist}{define}
	\item[$F :=$] $F_i$

	\item[$B :=$] $\kk[\ZZ F]$-basis for $(0:_N P_F)[\ZZ F]$, as
		in Algorithm~\ref{a:B}
	\end{routinelist}

	\routine{while}\ $y \in B$ and $g \in G$ \procedure{do}

	\begin{algsublist}

		\routine{if} $(0:_{\<g\>} P_F) \neq 0$ in some degree
			$a_{yg} \equiv \deg(y) \mod{\ZZ F}$

		\begin{algsublist}

			\begin{routinelist}{then}
			\item[$\lambda_{yg} :=$] scalar coefficient of
				$y$ on $\xx^{a_{yg}-\deg(g)} \cdot g$,
				as in Algorithm~\ref{a:B'}
			\end{routinelist}

			\routine{else} $\lambda_{yg} := 0$

		\end{algsublist}

		\routine{end} \procedure{if-then-else}

	\end{algsublist}

	\routine{end} \procedure{while-do}

	\begin{routinelist}{redefine}
	\item[$\lambda_g :=$] concatenation of the two vectors
		$\lambda_g$ and $(\lambda_{yg})_{y\in B}$, for $g\in G$
	\item[$\WW :=$] $\WW \oplus (\#B$ copies of $F,$
		$\ZZ^d$-degrees of vectors in $B)$

	\item[$N :=$] $M/\GF M$

	\item[$i :=$] $i+1$
	\end{routinelist}

\end{algsublist}

\routine{end} \procedure{while-do}

\routine{output} $\WW$ along with $\Lambda = \{\lambda_g\}_{g \in G}$,
	where $\lambda_g$ is in degree $\deg(g)$
\end{alglist}

\begin{prop} \label{p:irred}
Algorithm~\ref{a:irred} outputs an effective irreducible hull of~$M$,
using generators and relations for~$M$ as input.
\end{prop}
\begin{proof}
We must show that the homomorphism $M \to \WW$
determined by $G$ and $\Lambda$ is well-defined and injective.  More
precisely: monomial combinations $z$ of the generators of $M$ are zero
if and only if the corresponding monomial combinations $z_\lambda$ of
the $\lambda_g$ are zero in~$\WW$; here, $\lambda_g$ represents not a
data structure but an element of $\WW$.

The combination $z$ is nonzero in~$M$ if and only if the submodule
$\<z\> \subseteq M$ generated by $z$ has an associated prime.  The
associated prime is $F := F_i$ if and only if the image of $\<z\>$ in
the succesive quotient $N = M/\GF{}_{_{\!i-1}} M$ intersects $(0:_N
P_F)$ nontrivially (this in particular imlplies that $(0:_N P_F)$ is
nonzero, so $P_F$ is associated to~$M$).  This nontriviality of $\<z\>
\cap (0:_N P_F)$ is equivalent to having at least one of the terms
monomial$\cdot g$ appearing in~$z$ be nonzero in the same $(0:_N
P_F)$, because $B$ is a basis for $(0:_N P_F)[\ZZ F]$.  Finally,
monomial$\cdot g$ is nonzero precisely when the corresponding element
monomial$\cdot \lambda_g$ has nonzero coefficient in the appropriate
summand of~$\WW$.
\end{proof}

\begin{remark}\rm \label{r:avoid}
Some alterations to Algorithm~\ref{a:irred} may improve its running
time.
\begin{numbered}
\item
It is possible to avoid taking the successive quotients $N/\GF M$ at
the \procedure{redefine} step.  These quotients are designed to make
Lemmas~\ref{l:B} and~\ref{l:B'} apply, as well as to make $N$
successively simpler.  However, the cost of taking these quotients may
not be worth it, since the final sentence of Lemma~\ref{l:B} holds
even if $F$ doesn't have minimal dimension (so $(0:_M P_F)$ doesn't
include into its localization along~$F$).  In fact, both of
Algorithms~\ref{a:B} and~\ref{a:B'} still work in this more general
setting.

\item
Of the faces on the list $(F_1,\ldots,F_s)$, only those associated
to~$M$ need to be tested.
If desired, these faces can be detected using homological methods.

\item
Instead of computing and working with $(0:_M P_F)$ for each face
separately, one could work with the modules $(0:_M I_c)$ for each~$c$,
where $I_c$ is the intersection of all primes $P_F$ for faces $F$ of
dimension~$c$.
\end{numbered}
\end{remark}

\subsection{Generators and relations from irreducible hulls: saturated case}

\label{sub:satgenrel}

In this subsection we assume that $Q$ is saturated.  Our goal is to
compute relations on the generators for~$M$ that come as part of an
effective irreducible hull $M \into \WW$.  As we shall see in the
proof of Proposition~\ref{p:relations}, the computation essentially
reduces to the case where $M = \WW$ is an {\em indecomposable}
irreducible sum~$\WW$, so we are to determine the kernel of the
surjection $\kk[Q] \onto \WW$.
More explicitly, given a face $F$ and a degree $a \in Q$, we must find
generators of
\begin{eqnarray} \label{eq:W}
  W &:=& \kk\{Q \minus (a+F-Q)\}
\end{eqnarray}
as an ideal in $\kk[Q]$.

Since $Q$ is saturated, there is a unique minimal set of oriented
hyperplanes inside~$\ZZ^d$ whose closed positive half-spaces in
$\ZZ^d$ have intersection equal to~$Q$.  The map sending $H \mapsto H
\cap Q$ gives a bijection from these hyperplanes to the facets of~$Q$.
Denote by~$H_+$ the closed positive half-space determined by an
oriented hyperplane~$H$,
and by $H_+^\circ$ the {\em open} positive half-space.  Thus
$H_+^\circ$ is the complement of
$-H_+$ but can also be characterized as the lattice distance~$1$
translate of $H_+$ in the positive direction.

\begin{lemma} \label{l:hyperplanes}
Given any face $F$ of~$Q$ and any element $a \in Q$,
\begin{eqnarray*} \label{eq:H}
  Q \minus (a+F-Q) &=& \bigcup_{H \supseteq F} (a + H_+^\circ) \cap Q.
\end{eqnarray*}
\end{lemma}
\begin{proof}
We have
$a + F-Q = \bigcap_{H \supseteq F} a-H_+$ because $Q$ is saturated
(recall $F-Q = -(Q + \ZZ F)$).  Thus $\ZZ^d \minus (a + F-Q) =
\bigcup_{H \supseteq F} a + H_+^\circ$.  Now intersect with~$Q$.
\end{proof}

Lemma~\ref{l:hyperplanes} reduces the computation of generators for
$W$ as in~(\ref{eq:W}) to the case where $F$ is itself a facet, at
least when $Q$ is saturated.  The next algorithm and two lemmas cover
this case by producing some rational polytopes whose integer points do
the job.
For notation, $\RR_+ F$ denotes the real cone generated by $F$ in
$\RR^d = \RR \otimes \ZZ^d$, and $\RR H$ denotes the real span of a
hyperplane~$H$.  Also, by a \bem{$Q$-set} we mean a subset of~$\ZZ^d$
closed under addition by elements of~$Q$.  A~set~$G$ of vectors
in~$\ZZ^d$ \bem{generates} a $Q$-set $T$ if $T = G + Q$.

\begin{lemma} \label{l:magic}
Let $G_Q$ be the zonotope that is the Minkowski sum of all primitive
integer vectors along rays of~$Q$.  Then, for all $\alpha \in \RR^d$,
the lattice points in $\alpha + G_Q$ generate $(\alpha + \RR_+ Q) \cap
\ZZ^d$ as a $Q$-set.
\end{lemma}
\begin{proof}
Let $\beta$ be a lattice point in $\alpha + \RR_+Q$.  If there is no
primitive integer vector~$\rho$ along a ray of~$Q$ such that $\beta -
\rho$ still lies in $\alpha + \RR_+Q$, then $\beta \in \alpha + G_Q$.
\end{proof}

\begin{alg}\rm \label{a:polytopes}
\end{alg}
\begin{alglist}

\begin{routinelist}{input}
\item[$Q :=$] a saturated semigroup

\item[$H :=$] one of the hyperplanes bounding $Q$

\item[$a\,\,\in$] $Q$
\end{routinelist}

\routine{output} \ finite set $B \subset Q$ such that the ideal
$\<\xx^b \mid b \in B\>$ equals $\kk\{(a+H_+^\circ) \cap Q\}$

\begin{routinelist}{define}
\item[$G := $] the polytope $G_Q$ in Lemma~\ref{l:magic}

\item[$F := $] $H \cap Q$, a facet of~$Q$

\item[$\Delta :=$] the set of faces of $Q$ intersecting $F$ only at
$\0 \in Q$

\end{routinelist}

\begin{routinelist}{initialize}
\item[$B := $] $\{\}$, the empty subset of $Q$
\end{routinelist}

\routine{while} $\ D \in \Delta$ \procedure{do}

\begin{algsublist}

	\routine{define} $B_D :=$ lattice points in Minkowski sum
	$\bigl((a + \RR H) \cap \RR_+ D\bigr) + G$

	\routine{redefine} $B := B \cup B_D$

	\routine{next} $D$

\end{algsublist}

\routine{end} \procedure{while-do}

\routine{output} $B$
\end{alglist}

\begin{lemma} \label{l:polytopes}
Algorithm~\ref{a:polytopes} computes generators for the ideal
$\kk\{(a+H_+^\circ) \cap Q\}$.
\end{lemma}
\begin{proof}
Suppose $b \in (a + H_+^\circ) \cap Q$.  The intersection $(b + \RR H)
\cap \RR_+ Q$ is a polyhedron whose bounded faces are precisely the
polytopes $(b + \RR H) \cap \RR_+ D$ for $D \in \Delta$, and whose
recession cone is $\RR_+ F$.  Therefore $b \in b' + \RR_+ F$ for some
real vector $b' \in (b + \RR H) \cap \RR_+ D$ and some face $D \in
\Delta$.  Moreover, $b'$ lies in $b'' + \RR_+ D$ for some real vector
$b'' \in (a + \RR H) \cap \RR_+ D$.  Consequently, $b$ lies in $b'' +
\RR_+(D + F)$, and therefore in $b'' + \RR_+ Q$.  Now $\xx^b$ lies in
the $\kk[Q]$-module generated by $\kk\{B_D\}$, by definition of~$G$.
\end{proof}

\begin{remark}\rm \label{r:polytopes}
Some alterations to Algorithm~\ref{a:polytopes} may improve its
running time.
\begin{numbered}
\item
Instead of computing just one polytope $G = G_Q$ and Minkowski summing it
to define every~$B_D$, we could
define $B_D$ with $G_{D+F}$ in place of $G$, for each face $D \in
\Delta$.  This might reduce the number of lattice points in $B$
dramatically, but would require more computations as in
Lemma~\ref{l:magic}.

\item
Restricting to the maximal elements in $\Delta$ will speed things up.
\end{numbered}
\end{remark}

Let us summarize the above algorithm and three lemmas.  (See
Section~\ref{sec:injalg}
for issues concerning the output of the algorithm in the following
proposition, and post-processing for the purpose of reducing its
complexity.)

\begin{prop} \label{p:saturated}
Generators of the irreducible ideal\/ $W = \ker(\kk[Q] \onto \WW)$ are
algorithmically computable using as input an indecomposable effective
irreducible sum~$\WW$ over a normal semigroup ring~$\kk[Q]$.
\end{prop}
\begin{proof}
Apply Algorithm~\ref{a:polytopes} to each of the sets $(a + H_+^\circ)
\cap Q$ in Lemma~\ref{l:hyperplanes}.
\end{proof}

\subsection{Generators and relations from irreducible hulls: unsaturated case}

\label{sub:unsatgenrel}

Now we return to the general case, where~$Q$ need not be saturated,
and denote by~$Q^\sat$ the saturation of~$Q$.  The basic idea for
computing generators of irreducible ideals in~$\kk[Q]$ is to intersect
(as $\kk[Q]$-modules) the submodule $\kk[Q] \subset \kk[Q^\sat]$ with
the ideal $W \subseteq \kk[Q^\sat]$ output in the saturated case,
Proposition~\ref{p:saturated}.  Then it remains to find the
appropriate $F$-primary component of~$W$ as a $\kk[Q]$-module, where
$F$ is the unique face of dimension $\dim(\WW)$ associated to $\WW$
(as a $\kk[Q]$-module).

Every module in Algorithm~\ref{a:unsat} is to be considered as a
$\kk[Q]$-module---even those generated as $\kk[Q^\sat]$-modules.  Thus
$F$ is always a face of~$Q$, and we consider $F-Q$ as opposed to
$F-Q^\sat$.  Note, however, that $\kk\{F-Q^\sat\}$ does equal the
corresponding injective over $\kk[Q^\sat]$, even though $F$ is a face
of~$Q$; subtracting $Q^\sat$ automatically saturates~$F$.

\begin{alg} \label{a:unsat}
\end{alg}
\begin{alglist}

\begin{routinelist}{input}
\item[$Q :=$] a semigroup, not necessarily saturated

\item[$F :=$] a face of $Q$

\item[$a\,\,\in$] $Q$
\end{routinelist}

\begin{routinelist}{output}
\item[$B \subset$] $Q$ such that $\<\xx^b \mid b \in B\>$ equals the
ideal $\kk\{Q \minus (a+F-Q)\}$ in $\kk[Q]$
\end{routinelist}

\begin{routinelist}{define}
\item[$\ol V :=$] $\kk\{Q^\sat \minus (a + F - Q^\sat)\}$, an
indecomposable irreducible over~$\kk[Q^\sat]$

\item[$V := $] the kernel of $\kk[Q^\sat] \onto \ol V$ output by
Proposition~\ref{p:saturated}

\item[$W := $] $V \cap \kk[Q]$, the intersection taken inside
	$\kk[Q^\sat]$

\item[$I :=$] $\bigcap\,\{P_D \mid D$ is a facet of $F\}$, an ideal in
	$\kk[Q]$
\end{routinelist}

\begin{routinelist}{initialize}
\item[$B := $] degrees of the elements generating $W$

\item[$\WW := $] $\kk[Q]/W$
\end{routinelist}


\routine{while} $(0:_\WW P_F)$ has a generator in some degree
$\not\equiv a \mod{\ZZ F}$ \procedure{do}

\begin{algsublist}

	\begin{routinelist}{define}
	\item[$G := $] generators for $(0:_\WW P_F)$ that lie in
		degrees	$\not\equiv a \mod{\ZZ F}$
	\end{routinelist}

	\begin{routinelist}{redefine}
	\item[$B := $] $B \cup \mbox{}$degrees of the elements in $G$

	\item[$\WW :=$] $\WW/G$
	\end{routinelist}

	\begin{routinelist}{define}
	\item[$G' := $] generators for $\GI \WW$
	\end{routinelist}

	\begin{routinelist}{redefine}
	\item[$B := $] $B \cup \mbox{}$degrees of the elements in $G'$

	\item[$\WW :=$] $\WW/G'$
	\end{routinelist}

\end{algsublist}

\routine{end} \procedure{while-do}

\routine{output} $B$
\end{alglist}



\begin{prop} \label{p:W}
Algorithm~\ref{a:unsat} outputs generating degrees for $\kk\{Q \minus
(a+F-Q)\}$.
\end{prop}
\begin{proof}
The module $\WW$ gets initialized as a quotient of $\kk[Q^\sat]$ with
dimension~$\dim(F)$ as a $\kk[Q]$-module.  This much holds by the
saturated version Proposition~\ref{p:saturated} applied to $\ol V$,
and the preservation of dimension
\cite[Proposition~9.2]{Eis} for the module-finite ring extension
$\kk[Q] \subseteq \kk[Q^\sat]$ applied to~$V$.  One part of the output
is clear: the set $\<\xx^b \mid b \in B\>$ generates the kernel of the
map $\kk[Q] \onto \WW$ at every stage in the algorithm.  The question
is whether $\WW$ is the claimed indecomposable irreducible sum.

In the first \procedure{redefine} step,
the annihilator of $\xx^a \in \WW$ remains~$P_F$.  Indeed, any element
killed by~$P_F$ that generates a submodule containing a nonzero
element in degree~$a$ must itself have degree congruent to $a \mod{\ZZ
F}$.  The second \procedure{redefine} step only kills elements with
annihilators strictly larger than that of $\xx^a$; such elements can't
generate submodules containing~$\xx^a$.  Therefore, $\WW$ has only one
associated prime~$P_F$ after each loop of \procedure{while-do}, by
dimension considerations.

When the loop terminates, the localization $(0:_\WW P_F)[\ZZ F]$ along
$F$ is indecomposable, being isomorphic to $\kk\{a + \ZZ F\}$.  It
follows that the kernel of the surjection $\kk[Q] \onto \WW$ is an
irreducible ideal \cite[Proposition~3.1.7]{Vas}.  We are done by
Lemma~\ref{l:irrideal}, because $\kk\{a+F-Q\}_Q$ is the only
indecomposable irreducible sum for which the annihilator of $\xx^a$
is~$P_F$.
\end{proof}

\begin{remark}\rm
Some alterations to Algorithm~\ref{a:unsat} may improve its
efficiency.
\begin{numbered}
\item
The step $\WW := \WW/\GI \WW$ need not occur until the very last step
before \procedure{output}.  Its current placement is designed to speed
the computation by simplifying $\WW$ in each loop, but the cost of
taking the colon may not make up for it.  Instead, the end
of the algorithm can be replaced by:
\medskip

\begin{alglist}
\routine{while} $(0 :_\WW P_F)$ has rank strictly larger than $1$
over $\kk[F]$ \procedure{do}

\begin{algsublist}

	\begin{routinelist}{define}
	\item[$G :=$] generators for $(0:_\WW P_F)$ lying in degrees
		$\not\equiv a \mod{\ZZ F}$
	\end{routinelist}

	\begin{routinelist}{redefine}
	\item[$B :=$] $B \cup\mbox{}$degrees of the elements in $G$

	\item[$\WW :=$] $\WW/G$
	\end{routinelist}

\end{algsublist}

\routine{end} \procedure{while-do}

\begin{routinelist}{redefine}
\item[$B :=$] $B \cup\mbox{}$degrees of the generators of $\GI \WW$
\end{routinelist}

\routine{output} $\:B$
\end{alglist}
\medskip

\item
As in Remark~\ref{r:avoid}, it is not necessary to compute all of
$(0:_\WW P_F)$ in the \procedure{while-do} loop.  It suffices instead
to let $G$ be a basis for $(0:_\WW P_F)[\ZZ F]$.  This remark also
holds for the reworked \procedure{while-do} loop in the previous item.

\item
The set $B$ can become rather redundant.  Since the machine will have
to keep a presentation of $\WW$ in memory, the algorithm could
simply spit out the relations defining $\WW$ as a $\kk[Q]$-module at
the very end, without keeping track of $B$ at all.
\end{numbered}
\end{remark}

\begin{excise}{%
  How to compute generators for $Q^\sat$ as a $Q$-set:
  
  Pick any element $h \in Q^\sat$, and consider the semigroup
  generated by $h$ and~$Q$.  If $nh \in Q$ then $\0, h, 2h, \ldots,
  (n-1)h$ generate $\<Q,h\>$ as a $Q$-set because every element in
  $\<Q,h\>$ has the form $ih + a$ for some $i < n$ and $a \in Q$.  In
  terms of algebras, $\kk[Q,h] = \kk[Q][x_h]/$relations, and all the
  relations will have the form $\xx^a x_h = \xx^b$, where $a,b \in Q$
  and $a+h = b$ (proof: the relations obviously hold in $\kk[Q,h]$,
  and the induced map $\kk[Q][x_h]/$relations $\to \kk[Q,h]$ of
  $\kk[Q]$-modules is an isomorphism of vector spaces).
  
  If worse comes to worst, do the following now: each Hilbert basis
  element $h$ of $Q^\sat$ satisfies $n_h h \in Q$ for some $n_h \in
  \NN$ (easy to calculate: look for relation of the form $x_h^{n_h} -
  \xx^{\rm other\ stuff}$ in the toric ideal for $\<Q,h\>$ using any
  term order in which $h$ is really heavy, so that $x_h^{n_h}$ is in
  the initial ideal, or in which $\xx^{\rm other\ stuff}$ is so light,
  that it will never get chosen).  The set of all monomials $\prod_h
  x_h^{i_h}$ in which $i_h < n_h$ for each~$h$ generate $\kk[Q^\sat]$
  as a $\kk[Q]$-module.  The only remaining question is how to get the
  relations satisfied by these monomials.
  
  Start by choosing a subset of the monomials $m$ in the $x_h$
  variables in such a way that each monomial maps to a different
  element in $\kk[Q^\sat]$.  Take the direct sum of $\kk[Q]$-modules
  $\kk[Q,m]$ over all these monomials~$m$.  Mod out this direct sum by
  a couple of things.  First, kill all diagonals---that is, identify
  all of the copies of $\kk[Q]$ inside these $\kk[Q,m]$.  The result
  has $\kk$-vector space dimension~$1$ in all $Q$-graded degrees.  Now
  find the toric ideal for $Q^\sat$ using variables $x_h$ for Hilbert
  basis elements of $Q^\sat$ and other variables for~$Q$.  Mod out by
  all relations in this toric ideal that happen to lie in the degrees
  of the monomials~$m$, treating any product of the $x_h$ variables as
  the (now unique, because of the first sentence in this paragraph)
  monomial the multiply to.  That is, we get relations of the form
  $\xx^a m = \xx^b m'$ whenever $a + \deg(m) = b + \deg(m')$ in
  $Q^\sat$ and $a,b \in Q$.  The resulting quotient has vector space
  dimension~$1$ in every $Q^\sat$-graded degree, by construction.
  Since these relations all obviously hold in $\kk[Q^\sat]$, there is
  a natural induced map from the resulting quotient to $\kk[Q^\sat]$,
  and we have just seen that the map is an isomorphism.
  
  Note that, once we know {\em any} $\kk[Q^\sat]$-module $M$ as a
  $\kk[Q]$-module, we know how to express every
  $\kk[Q^\sat]$-submodule $M' \subseteq M$ as a $\kk[Q]$-module.
  Indeed, it suffices to give a generating set for $M'$.  But this is
  easy: take the generating set for $M'$ as a $\kk[Q^\sat]$-module
  (assume this has been given), and then simply multiply these by the
  generators of $\kk[Q^\sat]$ as a $\kk[Q]$-module.  This
  multiplication takes place over $\kk[Q^\sat]$, but the elements thus
  produced can be just as easily considered over $\kk[Q]$.  This is
  how we make any $\kk[Q^\sat]$-ideal of $\kk[Q^\sat]$ into a module
  over $\kk[Q]$.
}\end{excise}%
%

\section{Computing injective resolutions}

\label{sec:injres}

In this section the semigroup~$Q$ is not required to be saturated.
Our goal is the main result (Theorem~\ref{t:injres}) in the first half
of the paper: an algorithm to compute injective resolutions of
finitely generated modules over~$\kk[Q]$, in the $\ZZ^d$-graded
setting.  That is, given generators and relations for a finitely
generated $\ZZ^d$-graded module~$M$, we will compute an exact sequence
$0 \to M \to J^0 \to J^1 \to \cdots$\ in which $J^i$ is a
$\ZZ^d$-graded injective module for each~$i$.  Of course, we shall
only say how to calculate up to some specified cohomological degree,
as injective resolutions usually do not terminate.  This will not pose
a problem for our subsequent computation in Section~\ref{sec:injalg}
of local cohomology, which
vanishes past cohomological degree~\mbox{$d+1$} anyway.

The upshot is to reduce the computation of injective resolutions to
finding irreducible hulls of finitely generated $Q$-graded modules and
computing their cokernels, which we have already done in
Section~\ref{sec:irred}.

The data structures we employ for $\ZZ^d$-graded injective resolutions
are the matrices we introduce in the next definition.

\begin{defn}\rm \label{d:injmat}
A \bem{monomial matrix} is a matrix of constants~$\lambda_{qp}$ along
with
\begin{numbered}
\item
a vector~$\alpha_q \in \ZZ^d$ and a face~$F_q \in Q$ for each row, and

\item
a vector~$\alpha_p \in \ZZ^d$ and a face~$F_p \in Q$ for each column
\end{numbered}
such that $\lambda_{qp} = 0$ unless $F_p \subseteq F_q$ and $\alpha_p
\in \alpha_q + F_q-Q$.
\end{defn}
These monomial matrices generalize those in \cite{Mil2}, which were
for $Q = \NN^d$.

To any monomial matrix we can associate a map $J \mapsto J'$ of
injective modules in the following manner.  Each row and column label
gives the data of an indecomposable injective; we think of the row
labels as giving summands of~$J$ and the column labels as giving
summands of~$J'$.  To give a map from~$J$ to~$J'$ is thus the same as
giving a matrix of maps from the row indecomposables to the column
indecomposables.  Such a map $\kk\{\alpha_q + F_q-Q\} \mapsto
\kk\{\alpha_p + F_p-Q\}$ is necessarily zero unless $F_p \subseteq
F_q$ and $\alpha_p \in \alpha_q + F_q-Q$.  In the latter case it is
determined by a single scalar~$\lambda_{qp}$.  Hence
$$
\begin{array}{@{\ }c@{\ \;}c@{\ \;}c@{\ \;}}
  &
  \injmatrix
	{\vdots & \vdots\:\\F_q & \alpha_{q\!}\\\vdots & \vdots\:\\}
	{\begin{array}{c}
		\cdots\   F_p\    \cdots\\
		\cdots\ \alpha_p\ \cdots\\[.5ex]
					\\
		\lambda_{qp}		\\
					\\
	 \end{array}}
	{\\\\\\}
\end{array}
$$
is a monomial matrix representing a map
$$
  \bigoplus_q \kk\{\alpha_q + F_q-Q\}\ \mapsto\ \bigoplus_p
  \kk\{\alpha_p + F_p-Q\}.
$$
The component $\kk\{\alpha_q + F_q-Q\} \mapsto \kk\{\alpha_p +
F_p-Q\}$ of this homomorphism takes $\xx^\alpha$ to
$\lambda_{qp}\xx^\alpha$ for all $\alpha \in \alpha_p + F_p-Q$, and is
zero elsewhere.

Note that in degree $\alpha$, the map $J_{\alpha} \mapsto J'_{\alpha}$
given by a monomial matrix is obtained by deleting the rows and
columns labeled by $\alpha_p, F_p$ such that $\alpha$ does not lie in
$\alpha_p + F_p-Q$.  (This corresponds to ignoring those summands of
$J$ and $J'$ not supported at $\alpha$.)  Ignoring the labels on what
remains gives us a matrix with entries in $\kk$, which defines the
$\kk$-vector space map $J_{\alpha} \mapsto J'_{\alpha}$.

Two monomial matrices represent the same map of injectives (with given
decompositions into direct sums of indecomposable injectives) if and
only if (i)~their scalar entries are equal, (ii)~the corresponding
faces~$F_r$ are equal, where $r = p,q$, and (iii)~the corresponding
vectors~$\alpha_r$ are congruent modulo~$\ZZ F_r$.

Rather than compute directly with cumbersome, infinitely generated
injectives, it is more convenient to approximate injective resolutions
using irreducible sums.

\begin{defn}\rm \label{d:irredres}
An \bem{irreducible resolution} of a $Q$-graded module $M$ is an exact
sequence $0 \to M \to \WW^0 \to \WW^1 \to \cdots\ $ in which each
$\WW^j$ is an irreducible sum.
\end{defn}

Irreducible resolutions are approximations to injective resolutions;
indeed, the $Q$-graded part of any injective resolution is an
irreducible resolution \cite[Theorem~2.4]{MilIrr}.  In particular,
monomial matrices just as well represent homomorphisms of irreducible
sums, as long as the degree labels $\alpha_q$ and~$\alpha_p$ all can
be chosen to lie in~$Q$.  The (apparent) advantage to irreducible
resolutions over injective resolutions is their finiteness.

\begin{cor} \label{c:irredres}
For any finitely generated $Q$-graded $\kk[Q]$-module\/~$M$,
Propositions~\ref{p:relations} and~\ref{p:irred} inductively compute a
minimal irreducible resolution $\WW^\spot$ of\/~$M$
\mbox{algorithmically}.
\end{cor}
\begin{proof}
Minimal irreducible resolutions have finite length (that is, they
vanish in all sufficiently high cohomological degrees) by
\cite[Theorem~2.4]{MilIrr}.  The computability therefore follows from
Propositions~\ref{p:relations} and~\ref{p:irred} by induction on the
highest cohomological degree required.
\end{proof}


The next result demonstrates the precise manner in which irreducible
resolutions approximate injective resolutions for computational
purposes.

\begin{prop} \label{p:Qgraded}
Let $M$ be a finitely generated module with minimal injective
resolution $J^\spot$ and minimal irreducible resolution $\WW^\spot$.
Suppose that every indecomposable summand in the first $n$
cohomological degrees of $J^\spot$ has nonzero $Q$-graded part.  Then
$M$ is $Q$-graded, and the data contained in the first $n$ stages of
$\WW^\spot$ constitute a finite data structure for the first $n$
cohomological degrees of~$J^\spot$.
\end{prop}
\begin{proof}
Every map in $J^\spot$ can be expressed using the finite data of a
monomial matrix, and this data can be read immediately off the maps
in~$\WW^\spot$.
\end{proof}

If we can algorithmically determine a $\ZZ^d$-graded shift of $M$ so
that the hypotheses of Proposition~\ref{p:Qgraded} are satisfied, then
we can compute the minimal injective resolution of~$M$ up to
cohomological degree~$n$.  This task requires a lemma, in which $\mm$
denotes the maximal ideal $P_{\{\0\}}$ generated by nonunit monomials
in~$\kk[Q]$.

\begin{lemma} \label{l:pinnacle}
Let $J^\spot$ be a minimal injective resolution of a finitely
generated module $M$, and $F$ a face of\/~$Q$.  If every
indecomposable summand of\/ $\Gm J^{j+d-\dim(F)}$ has nonzero
$Q$-graded part, then every indecomposable summand of~$J^j$ isomorphic
to a $\ZZ^d$-graded shift of\/ $\kk\{F-Q\}$ has nonzero $Q$-graded
part.
\end{lemma}
\begin{proof}
\cite[Proposition~3.5]{HelM}, in the special case of an affine
semigroup ring.
\end{proof}

Every indecomposable summand of $\Gm J^j$ is a shift $\kk\{\alpha-Q\}$
of~$\kk\{-Q\}$.  Such an indecomposable injective has nonzero
$Q$-graded part if and only if $\alpha \in Q$.
Our final lemma in this section describes the (standard) way to
calculate the shifts $\alpha$ appearing in $\Gm J^j$.  The number
$\mu^{j,\alpha}(M)$ of shifts $\kk\{\alpha-Q\}$ appearing as summands
in cohomological degree $j$ of the minimal injective resolution of~$M$
is called the $j^\th$ \bem{Bass number} of~$M$ in degree $\alpha$.

\begin{lemma} \label{l:bass}
Let $\FF_\spot$ be a free resolution of the residue field~$\kk$.  The
Bass number $\mu^{j,\alpha}(M)$ is effectively computable as the
$\kk$-vector space dimension of $H^j(\FF_\spot, M)_\alpha$.
\end{lemma}
\begin{proof}
This expression of Bass numbers as dimensions (over $\kk$) of Ext
modules is standard; see \cite[Chapter~3]{BH}.  The computability
follows because we can calculate free resolutions, homomorphisms, and
homology over~$\kk[Q]$.
\end{proof}

Now we come to our central result.  For notation, $M(-a)$ denotes the
$\ZZ^d$-graded shift of~$M$ up by~$a$, so that $M(-a)_b = M_{b-a}$.

\begin{thm} \label{t:injres}
Fix a finitely generated $\kk[Q]$-module $M$ and an integer~$i$.
There is an algorithmically computable $a \in Q$ for which
Propositions~\ref{p:relations} and~\ref{p:irred} inductively compute
the minimal injective resolution of~$M(-a)$ through cohomological
degree~$i+1$.
\end{thm}
\begin{proof}
After using Lemma~\ref{l:bass} to compute the Bass numbers of $M$ up
to
cohomological degree $i+1+\dim(M)$, choose $a$ so that the
corresponding Bass numbers of $M(-a)$ have $\ZZ^d$-graded degrees
lying in~$Q$.  At this point, $M(-a)$ satisfies the hypotheses of
Proposition~\ref{p:Qgraded} with $n=i+1$, by Lemma~\ref{l:pinnacle}.
Now apply Corollary~\ref{c:irredres}.
\end{proof}

\section{Sector partitions from injectives}

\label{sec:decomp}


We turn now to sector partitions, for which we assume henceforth that
the affine semigroup $Q$ is saturated.  As a prerequisite to producing
sector partitions of local cohomology modules, we demonstrate in this
section that injective modules admit sector partitions, as does the
homology of any complex of injective modules.

\begin{prop} \label{p:sectors}
Suppose~$J = \bigoplus_{i=1}^r J_i$ is an injective module decomposed
into summands \mbox{$J_i = \kk\{\alpha_i + F_i-Q\}$}.  For each subset
$A \subseteq \{1, \dots, r\}$ define~$S_{\!A}$ to be the set
\begin{eqnarray*}
  S_{\!A} &=& \{\alpha \in \ZZ^d \mid (J_i)_\alpha \cong \kk \hbox{
  for } i \in A\}
\end{eqnarray*}
of all degrees in\/~$\ZZ^d$ such that the summands of~$J$ nonzero in
that degree are precisely those indexed by~$A$.
The sets~$S_{\!A}$ canonically determine a sector partition $\SS(J)
\vdash J$.
\end{prop}
\begin{proof}
For each $\alpha \in \ZZ^d$, either $(J_i)_\alpha = \{0\}$ or
$(J_i)_\alpha = \kk\cdot\xx^{\alpha}$.  Therefore $\SS(J)$ is indeed a
partition of~$\ZZ^d$.  Now we must show that $S_{\!A}$ is a finite
union of polyhedra as in part~1 of Definition~\ref{d:sector}.  The set
$\alpha + F-Q$ of degrees is the set of lattice points in a polyhedron
of the desired form because the half-spaces whose intersection is
$\alpha + F-Q$ are bounded by hyperplanes parallel to facets of~$Q$,
by definition.  These hyperplanes divide~$\ZZ^d$ into finitely many
disjoint regions (place the lattice points lying on each hyperplane in
the region on the positive side of that hyperplane), each of which
consists of the lattice points in a polyhedron of the desired form.
Thus the complement $\ZZ^d \minus (\alpha + F-Q)$ is the required kind
of finite union.  We conclude that $S_{\!A}$ is a finite union of
regions, each of which is an intersection of $r$ polyhedral
regions---one from each of the summands~$J_i$.

For each index set~$A$ such that $S_{\!A}$ is nonempty, define
$J_{S_{\!A}} \subseteq \kk^r$ to be the subspace spanned by the basis
vectors~$e_i$ such that $i \in A$.  Then for each degree~$\alpha$
in~$S_{\!A}$, the map $J_\alpha \to J_{S_{\!A}}$ required by part~2 of
Definition~\ref{d:sector} can be taken to equal the zero map on
$(J_i)_\alpha$ for $i$ not in~$A$, and the map sending $\xx^\alpha$
to~$e_i$ on $(J_i)_\alpha$ for $i$ in~$A$.

To define the maps $\xx^{S_B - S_{\!A}}$ for index sets $A$ and~$B$
such that $S_B - S_{\!A}$ is nonempty, as in part~3 of
Definition~\ref{d:sector}, it suffices to define the image of~$e_i$
for each~$i$ in~$A$.  We take $\xx^{S_B - S_{\!A}}(e_i) = e_i$ if $i$
is in~$B$, and $\xx^{S_B - S_{\!A}}(e_i) = 0$ otherwise.
Commutativity of the required diagram follows from the definition of
the module structure on $\kk\{\alpha_i + F_i-Q\}$.  Specifically, for
$\alpha \in S_{\!A}$ and $\beta \in S_B$ with $\beta - \alpha \in Q$,
multiplication by $\xx^{\beta - \alpha}$ takes $\xx^{\alpha}$ to
$\xx^{\beta}$ in $J_i$ for $i \in B$, and takes $\xx^{\alpha}$ to zero
in $J_i$ for $i$ outside~$B$.
\end{proof}

The sector partition in Proposition~\ref{p:sectors} descends to the
cohomology~$H$ of any complex of injectives, via monomial matrices.
The forthcoming sector partition
of~$H$ is really determined {\em canonically}\/ by~$J^\spot$ (without
its direct sum decomposition), even though the way we present things
here makes it look like bases must be chosen.  We chose this route
because bases are good for computation, while uniqueness is
immaterial.

\begin{thm} \label{t:decomp}
If $H$ is a module that can be expressed as the (middle) homology of a
complex $J^\spot: J' \to J \to J''$ in which all three modules are
injective, or all three modules are flat, then there is a
sector partition $\SS(J^\spot) \vdash H$ determined
by~$J^\spot$.
\end{thm}
\begin{proof}
Choose direct sum decompositions to write
$$
\begin{array}{r@{\ }c@{\ }l@{\quad\ }r@{\ }c@{\ }lcr@{\ }c@{\ }l}
  J'   &=& \dis\bigoplus_{i=1}^{r'} J'_i,
&
  J    &=& \dis\bigoplus_{i=1}^r J_i,
&\hbox{and}&
  J''\!&=& \dis\bigoplus_{i=1}^{r''} J''_i.
\end{array}
$$
Let $\Phi$ and $\Psi$ be the monomial matrices representing the maps
$J' \mapsto J$ and $J \mapsto J''$, respectively.  The sectors in the
sector partition $\SS(J' \oplus J \oplus J'') \vdash J' \oplus J
\oplus J''$ are indexed by triples $(A',A,A'')$ of subsets of
$\{1,\ldots,r'\}, \{1,\ldots,r\}, \{1,\dots,r''\}$, respectively, and
automatically satisfy the polyhedrality condition in part~1 of
Definition~\ref{d:sector} by Proposition~\ref{p:sectors}.
We take $\SS(J^\spot)$ to partition $\ZZ^d$ into these sectors.

For each triple $(A',A,A'')$ we have maps $\Phi_{A'}^A:
J'{}_{\!\!\!S_{\!A'}} \rightarrow J_{S_{\!A}}$ and $\Psi_{A}^{A''}:
J_{S_{\!A}} \to J''{}_{\!\!\!\!S_{\!A''}}$ whose monomial matrices are
defined by deleting: row~$i'$ of~$\Phi$ for $i'$ not in~$A'$;
column~$i$ of~$\Phi$ and row~$i$ of~$\Psi$ for $i$ not in~$A$; and
column~$i''$ of~$\Psi$ for $i''$ not in~$A''$.
Let
\begin{eqnarray} \label{eq:A}
  H_{S_{\!A',A,A''}} &=& \ker(\Psi_A^{A''})/\im(\Phi_{A'}^A).
\end{eqnarray}
For any $\alpha$ in $S_{\!A',A,A''}$, we have a commutative diagram
\begin{equation} \label{eq:AA}
\begin{array}{c@{\ }c@{\ }c@{\ }c@{\ }c}
J'_{\alpha}&\longrightarrow&J_{\alpha}&\longrightarrow&J''_\alpha\\
\downarrow &               &\downarrow&               &\downarrow\\[-1.5ex]
J'_{S_{A'}}&\stackrel{\Phi_{A'}^A}\too& J_{S_{\!A}}
           &\stackrel{\Psi_A^{A''}}\too&J''_{S_{A''}} 
\end{array}
\end{equation}
that induces the required isomorphism $H_\alpha \cong
H_{S_{\!A',A,A''}}$.  It is routine to check that the maps
$H_{S_{\!A',A,A''}} \to H_{S_{B',B,B''}}$ induced from the
corresponding maps on $J'_{A'}, J_A,$ and~$J''_{A''}$ commute with
this isomorphism.
\end{proof}

Once we have Theorem~\ref{t:decomp}, the only step remaining to prove
Theorem~\ref{t:sector} is to exhibit $H^i_I(M)$ as the homology of a
complex of injectives.

\begin{remark}\rm
The results in this section hold just as well for \bem{flat} objects
of~$\MM$, which are Matlis dual to injective objects and hence
isomorphic to finite direct sums of modules of the form $\kk\{\alpha +
F + Q\}$ for some $\alpha$ in~$\ZZ^d$ and some face~$F$ of~$Q$
\cite[Chapter~11]{cca}.
For the proofs, simply apply Matlis duality to the results for
injectives.
\end{remark}

\section{Computing sector partitions}

\label{sec:secposet}


Again letting $Q$ be a saturated affine semigroup, the next task is
actually computing the finitely many polyhedra whose lattice points
comprise the sectors in the sector partition $\SS(J) \vdash J$ of an
injective module.  That is, we need to make
Proposition~\ref{p:sectors} and its proof into an algorithm.

Since $Q$ is saturated, there are unique primitive integer linear
functionals $\tau_1,\ldots,\tau_n$ taking $\ZZ^d \to \ZZ$, one for
each facet of~$Q$, such that $Q = \bigcap_{i=1}^n \{\tau_i \geq 0\}$
is the set of lattice points in the intersection of their positive
half-spaces.  The degrees on which indecomposable injectives are
supported can be expressed in terms of these linear functionals, via
the following identity:
\begin{eqnarray} \label{eq:cutout}
  \alpha + F-Q &=& \{\beta \in \ZZ^d \mid \tau_i(\beta) \leq
  \tau_i(\alpha) \hbox{ whenever } \tau_i(F) = 0\}.
\end{eqnarray}
In other words, $F-Q$ is the intersection of the {\em negative}
half-spaces for those functionals~$\tau_i$ vanishing on~$F$, and
$\alpha + F-Q$ is simply a translate.  By convention,
we use the notation $\tau_i(\beta) \leq \infty$ to mean that there is
no restriction on the value of $\tau_i(\beta)$.  This allows a
notation $\tau_F(\alpha) \in (\ZZ\cup\infty)^n$ for the vector whose
$i^\th$ coordinate satisfies
\begin{eqnarray*}
   \tau_F(\alpha)_i &=& \left\{\begin{array}{@{}ll}
			  \tau_i(\alpha) & \hbox{if }\tau_i(F) = 0\\
			  \infty & \hbox{otherwise.}
			\end{array}\right.
\end{eqnarray*}
The point is that a vector $\beta \in \ZZ^d$ lies in $\alpha + F-Q$ if
and only if $\tau(\beta) \leq \tau_F(\alpha)$, where
\begin{eqnarray*}
  \tau(\beta) &=& \big(\tau_1(\beta), \ldots, \tau_n(\beta)\big)
\end{eqnarray*}
and the `$\leq$' symbol denotes componentwise comparison.  We shall
use the corresponding definitions of $\tau_F(\alpha)$ and
$\tau(\beta)$ for vectors $\alpha,\beta \in \RR^d = \RR \otimes
\ZZ^d$, so $\tau_F(\alpha) \in (\RR\cup\infty)^n$.

For the rest of this section, let
\begin{eqnarray} \label{eq:J}
  J &=& \bigoplus_{j=1}^r J^j,\quad\hbox{with}\quad J^j\ =\
  \kk\{\alpha_j + F_j-Q\},
\end{eqnarray}
be an injective module,
and define
\begin{eqnarray*}
  \tau^j &:=& \tau_{F_j}(\alpha_j) \quad\hbox{for }
  j = 1,\ldots,r.
\end{eqnarray*}
Thus for $i = 1,\ldots,n$ the vector $\tau^j$ has $i^\th$ coordinate
$\tau^j_i = \tau_{F_j}(\alpha_j)_i$, which equals either
$\tau_i(\alpha_j)$ or~$\infty$, depending on whether $\tau_i$ vanishes
on~$F_j$ or not.  Even without calculating the set~$\SS(J)$
algorithmically, the vectors~$\tau^j$ specify the map from~$\ZZ^d$
to~$\SS(J)$, by definition.  We record a precise version of this
statement in the next lemma.

\begin{lemma} \label{l:store}
A degree $\alpha \in \ZZ^d$ lies in $S_{\!A}$ if and only if \mbox{$A
= \big\{j\!\in\!\{1,\ldots,r\} \mid \tau(\alpha) \leq \tau^j\big\}$}.
\end{lemma}

It remains to ascertain which sets~$S_{\!A}$ of lattice points are
nonempty, and to determine the pairs $A,B$ for which we must compute a
map $\xx^{B-A} : J_A \to J_B$.  (The maps themselves, which are
canonical,
are constructed in the proof of Proposition~\ref{p:sectors}.)  For
each functional~$\tau_i$ there is a permutation~$w_i$ of
$\{1,\ldots,r\}$ satisfying $\tau_i^{w_i(1)} \leq \cdots \leq
\tau_i^{w_i(r)}$.  To simplify notation, we write $\tilde\tau_i^\ell$
instead of $\tau_i^{w_i(\ell)}$.  Also, set $\tilde\tau_i^0 = -\infty$
and $\tilde\tau_i^{r+1} = \infty$.

For fixed~$i$, the parallel affine hyperplanes $\{\tau_i =
\tilde\tau_i^\ell\}_{\ell=1}^r$ divide $\ZZ^d$ into strips
\begin{eqnarray*}
  \{\beta \in \ZZ^d \mid \tilde\tau_i^\ell + 1 \leq \tau_i(\beta) \leq
  \tilde\tau_i^{\ell+1}\}
\end{eqnarray*}
for $\ell = 0,\ldots,r$.  At most $r+1$ of these strips are nonempty,
because some of the hyperplanes may coincide.  Also, the last few of the
$\tilde\tau_i^\ell$ will equal~$\infty$; we interpret any strip where
$\tau_i^\ell = \tau_i^{\ell+1} = \infty$ as empty, and ignore it.

\begin{prop} \label{p:strips}
Let $J$ be as in~(\ref{eq:J}).
For any fixed $\ell_1,\ldots,\ell_n \in \{0,\ldots,r\}$, the lattice
points in the polyhedron
\begin{eqnarray*}
  \Delta(\ell_1,\ldots,\ell_n) &:=& \bigcap_{i=1}^n\{\beta \in \RR^d
  \mid \tilde\tau_i^{\ell_i} + 1 \leq \tau_i(\beta) \leq
  \tilde\tau_i^{\ell_i+1}\}
\end{eqnarray*}
all lie inside a single sector in\/ $\SS(J)$.  The partition
of\/~$\ZZ^d$ by the polyhedra $\Delta(\ell_1,\ldots,\ell_n)$ refines
the partition of\/~$\ZZ^d$ by the sectors in~$\SS(J)$.
\end{prop}
\begin{proof}
This follows from the definitions and~(\ref{eq:cutout}), which uses
that $Q$ is saturated.
\end{proof}

Proposition~\ref{p:strips} makes way for an algorithm to compute the
set of sectors.

\begin{alg}\rm \label{a:secset}
\end{alg}
\begin{alglist}
\routine{input} $J = \bigoplus_{j=1}^r J^j$, an injective module over
	$\kk[Q]$, with $J^j = \kk\{\alpha_j + F_j-Q\}$
\routine{output} the set $\SS(J)$ of sectors, each expressed as a list
	of polyhedra that partition it
\begin{routinelist}{define}
\item[$\phi : \ZZ^d \to$] subsets of $\{1,\ldots,r\}$, as in
	Lemma~\ref{l:store}
\end{routinelist}
\begin{routinelist}{initialize}
\item[$\AA := \{\}$,] the empty collection of subsets of $\{1,\ldots,r\}$
\end{routinelist}
\routine{while} $\ell_1,\ldots,\ell_n \in \{0,\ldots,r\}$ \procedure{do}
\begin{algsublist}
	\routine{if} $\Delta(\ell_1,\ldots,\ell_n)\neq\nothing$
	\begin{algsublist}
 		\routine{then} \procedure{define} $A :=
			\phi(\Delta(\ell_1,\ldots,\ell_n))$
		\routine{else} \procedure{next} $(\ell_1,\ldots,\ell_n)$
	\end{algsublist}
	\routine{end} \procedure{if-then-else}
	\routine{if} $A \in \AA$
	\begin{algsublist}
		\begin{routinelist}{then}
			\routine{redefine} $S_{\!A} := S_{\!A} \cup
				\{\Delta(\ell_1,\ldots,\ell_n)\}$
		\end{routinelist}
		\begin{routinelist}{else}
			\routine{initialize} $S_{\!A} :=
				\{\Delta(\ell_1,\ldots,\ell_n)\}$
			\routine{redefine} $\AA := \AA\cup\{A\}$
		\end{routinelist}
	\end{algsublist}
	\routine{end} \procedure{if-then-else}
	\routine{next} $(\ell_1,\ldots,\ell_n)$
\end{algsublist}
\routine{end} \procedure{while-do}
\routine{output} $\{S_{\!A} \mid A \in \AA\}$
\end{alglist}

\medskip
\noindent
Note that $\phi$ is constant on $\Delta(\ell_1,\ldots,\ell_n)$ by
definition, and can easily be determined directly from the data
$(\ell_1,\ldots,\ell_n)$.

Next comes the determination of which maps $\xx^{B-A}$ need computing.
In the coming algorithm, we write $\Delta(\ell_1,\ldots,\ell_n) \leq
\Delta(\ell_1',\ldots,\ell_n')$ if $(\ell_1,\ldots,\ell_n) \leq
(\ell_1',\ldots,\ell_n')$ as vectors in $(\ZZ\cup\infty)^n$.  Such
notation is justified because $\Delta(\ell_1,\ldots,\ell_n) \not\leq
\Delta(\ell_1',\ldots,\ell_n')$ automatically implies that
$\Delta(\ell_1',\ldots,\ell'_n) - \Delta(\ell_1,\ldots,\ell_n)$ fails
to intersect~$Q$.

\begin{alg}\rm \label{a:secposet}
\end{alg}
\begin{alglist}
\routine{input} sectors $S_{\!A}$ and $S_B$ in $\SS(J)$ from the
	output of Algorithm~\ref{a:secset}
\routine{output} the truth value of:
	``there exist $\alpha \in S_{\!A}$ and $\beta \in S_B$ with
	$\beta-\alpha \in Q$''
\routine{initialize} $\mathit{val} := $ \procedure{false}
\routine{while} $(\Delta_A,\Delta_B) \in S_{\!A} \times S_B$
	\procedure{and} $\mathit{val} = $ \procedure{false},
	\procedure{do}
\begin{algsublist}
	\routine{if} $A \supseteq B$ \procedure{and} $\Delta_A \leq
		\Delta_B$
	\begin{algsublist}
		\routine{then} \procedure{define} $\Delta_B - \Delta_A:=$
			the Minkowski sum of $\Delta_B$ and $-\Delta_A$
		\routine{else} \procedure{next} $(\Delta_A,\Delta_B)$
	\end{algsublist}
	\routine{end} \procedure{if-then-else}
	\routine{if} $Q \cap (\Delta_B - \Delta_A) \neq \nothing$
	\begin{algsublist}
		\routine{then} \procedure{redefine} $\mathit{val}:=$
			\procedure{true}
		\routine{else} \procedure{next} $(\Delta_A,\Delta_B)$
	\end{algsublist}
	\routine{end} \procedure{if-then-else}
\end{algsublist}
\routine{end} \procedure{while-do}
\routine{output} $\mathit{val}$
\end{alglist}

The proof of correctness for Algorithm~\ref{a:secposet} is
straightforward from the definitions, except for the first
\procedure{if-then-else} procedure, which relies on
Lemma~\ref{l:poset}, below.  Note the non-necessity in
Algorithm~\ref{a:secposet} of actually finding a witness in $\Delta_B
- \Delta_A$ for $S_{\!A} \preceq S_B$; as we have seen in (\ref{eq:A})
and~(\ref{eq:AA}) from the proof of Theorem~\ref{t:decomp}, the
natural map on cohomology is induced by taking submatrices of the
monomial matrix, regardless of where the witnesses lie.

\begin{lemma} \label{l:poset}
If $\SS(J)$ is as in Proposition~\ref{p:sectors}, then $Q \cap (S_{B}
- S_{\!A}) \neq \nothing$ implies $A \supseteq B$.
\end{lemma}
\begin{proof}
If $a \in Q$ and $(J_i)_\alpha = 0$, then $(J_i)_{a+\alpha} = 0$, so the
set of summands nonzero in degree $a+\alpha$ can only be smaller.
\end{proof}

Unfortunately, Algorithm~\ref{a:secposet} is necessary, because
\hbox{$\Delta_B - \Delta_A \neq \nothing$} need not always hold when
$\Delta_A \leq \Delta_B$, as the example to come shortly demonstrates.
It does seem, however, that the offending pairs of polytopes are
usually ``small''.  For instance, we know of no examples where the
lattice points in either polytope affinely span~$\ZZ^d$.

\begin{example}\rm
Let $Q$ be the subsemigroup of $\NN^2$ generated by $(2,0)$, $(1,1)$,
and $(0,2)$.  Name the faces of $Q$ as $\0, X, Y, Q$, and set $E_F =
F-Q$.  Let
\begin{eqnarray*}
  J &=& \kk\{(0,0)+E_\0\}\,\oplus\,\kk\{(0,1)+E_X\}\,\oplus\,\kk\{(0,0)+E_Y\}\\
    & &\phantom{\kk\{(0,0)+E_\0\}}\oplus \kk\{(0,-1)+E_X\}\oplus
    \kk\{(-2,0)+E_Y\},
\end{eqnarray*}
with the summands labeled in order as $J_1,\ldots,J_5$.  Letting $X$ be
facet number~$1$ and $Y$ be facet number~$2$, the arrays $\tau_i^j$ and
$\tilde\tau_i^\ell$ look like
\begin{rcgraph}
  \left(\begin{array}{@{\,}c@{\,}}\tau_1^j\\\tau_2^j\end{array}\right)
  \ =\
  \left(
  \begin{array}{@{\;}ccccc@{\;}} 0 & 1 & \infty & -1 & \infty\\[.5ex]
		                 0 & \infty & 0 & \infty & -2
  \end{array}
  \right)
\qquad\hbox{and}\qquad
  \left(
  \begin{array}{@{\,}c@{\,}}\tilde\tau_1^\ell\\[.5ex]
			    \tilde\tau_2^\ell\end{array}
  \right)
  \ =\
  \left(
  \begin{array}{@{\;}ccccccc@{\;}}
		 -\infty & -1 & 0 & 1 & \infty & \infty & \infty\\[.5ex]
		 -\infty & -2 & 0 & 0 & \infty & \infty & \infty
  \end{array}
  \right)
\end{rcgraph}
The sectors $S_{\{1,2,3\}}$ and $S_{\{2,3\}}$ contain one polytope
each, and both of these polytopes contain exactly one lattice point.
Specifically, identifying the sector, the polytope, and the lattice
point, we have
$$
  S_{\{1,2,3\}} = \Delta(-1,-2) = (0,0)  \qquad\hbox{and}\qquad
  S_{\{2,3\}} = \Delta(0,-2) = (-1,1).
$$
Now $\Delta(-1,-2) \leq \Delta(0,-2)$, but subtracting the vector
in~$S_{\{1,2,3\}}$ from the one in~$S_{\{2,3\}}$ yields $(-1,1)$,
which does not lie in the semigroup~$Q$.%
\end{example}

\begin{remark}\rm \label{rk:holes}
The notion of sector partition ought to have a refinement that takes
into account the various kinds of failures of saturation for arbitrary
affine semigroup.  The resulting notion would produce sector
partitions for the cohomology of complexes of injectives over
nonnormal affine semigroup rings.  The failures of saturation fall
into two categories: the geometric kind, arising from polyhedral
``holes'' in the semigroup (as compared with its saturation), and the
arithmetic kind, arising from finite-index sublattices generated by
faces.  Even in the case where arithmetic failure is absent, however,
we do not know how to bound the sizes and shapes of the ``holes''
sufficiently to carry out an analysis such as the one producing the
algorithms above.
%
\end{remark}

\section{Computing local cohomology with monomial support}

\label{sec:injalg}

Still assuming that $Q$ is sturated, we have now finally developed
enough tools to prove the main theorem on local cohomology with
monomial support, namely Theorem~\ref{t:sector} from the Introduction.

\begin{proofof}{Theorem~\ref{t:sector}}
Take $i = d$ in Theorem~\ref{t:injres}, and let $J^\spot(-a)$ be the
minimal $\ZZ^d$-graded injective resolution computed there.  Then
$J^\spot$ is an algorithmically computed injective resolution of~$M$.
By definition, $H^i_I(M)$ is the middle cohomology of the complex
$\GI J^{i-1} \to \GI J^i \to \GI J^{i+1}$, where $\GI J^j$ is the
direct sum of all indecomposable summands of $J^j$ whose unique
associated prime contains~$I$.  Having now expressed $H^i_I(M)$ as the
cohomology of an effectively computed complex of injectives,
Theorem~\ref{t:decomp} says that $H^i_I(M)$ has a sector partition.
The set of sectors in part~1 of Definition~\ref{d:sector} is computed
by Algorithm~\ref{a:secset}.  The vector spaces in part~2 of
Definition~\ref{d:sector} are specified in (\ref{eq:A}) from the proof
of Theorem~\ref{t:decomp}, and naturally determine the maps in part~3
of Definition~\ref{d:sector}, given the computation in
Algorithm~\ref{a:secposet}.
\end{proofof}

%



Now we turn to issues of complexity.  There is little sense in
completing a formal complexity analysis of all of the algorithms
presented in this paper, as they involve Gr\"obner basis computation,
which is doubly-exponential from a worst-case perspective.  However,
it is worth mentioning where the complexity in our algorithms comes
from, up to a factor arising from the complexity of Gr\"obner basis
computation, since Gr\"obner basis computations are often more
efficient than expected.  The purpose of what follows, therefore, is
to assure the reader that our algorithms have not amplified the
faux-doubly-exponential complexity of Gr\"obner bases with some
``honest'' exponential complexity.

Let us assume that the dimension~$d$ is fixed, and analyze the
complexity of computing all of the local cohomology of a finitely
generated module~$M$ supported on a fixed monomial ideal~$I$ over a
normal semigroup ring~$\kk[Q]$.  This computation involves all of the
algorithms in the paper except the one in
Section~\ref{sub:unsatgenrel}.  (The complexity of
Algorithm~\ref{a:unsat} above and beyond Algorithm~\ref{a:irred} is
only about as bad as that of $\kk[Q^\sat]/\kk[Q]$ as a
$\kk[Q]$-module, anyway.)

In Algorithm~\ref{a:irred}, the only non-Gr\"obner contribution to the
running time comes from the number of basis elements constructed (see
Remark~\ref{r:avoid}.2, which can be used to ensure that we only check
faces of~$Q$ giving rise to basis elements).  This number is by
definition a Bass number of~$M$.  Thus, up to Gr\"obner basis
computation, Algorithm~\ref{a:irred} is only as complex as its output.

Next we consider
the algorithm in Proposition~\ref{p:saturated}.  The algorithm works
by taking the union (over a set of facets of~$Q$) of ideals output by
Algorithm~\ref{a:polytopes}.  The output presents the generators of
each such ideal as the lattice points in a union of polytopes having
the form $\bigl((a + \RR H) \cap \RR_+ D\bigr) + G$, where $D$ is a
face of~$Q$.
The computation of each such polytope is by standard techniques to
intersect polyhedra and take Minkowski sums with the fixed
zonotope~$G$.  Hence, up to factors coming from the number of facets
of~$Q$ and from standard procedures, we need only bound
\begin{numbered}
\item
the number of polytopes output by Algorithm~\ref{a:polytopes}, and

\item
the number of lattice points in each such polytope.
\end{numbered}
The former is polynomial in the number of facets of~$Q$ by
Remark~\ref{r:polytopes}.2.  The latter is polynomial in the input
vector $a \in Q$ by the piecewise polynomiality of the
lattice point enumeration function of $(a + \RR H) \cap \RR_+ D$ as a
function of~$a$ \cite{McMullen77}, along with the fact that $G$ is
fixed.  Actually computing the set of lattice points in each polytope
can be accomplished using the efficient algorithms of Barvinok and
Woods \cite{BarWoods}.

\begin{remark}\rm \label{r:BarWoods}
We need to do Gr\"obner basis computations with the irreducible
ideals~$W$ output by the algorithm in Proposition~\ref{p:saturated}.
This means that, for our purposes, the {\em short rational generating
functions}\/ output by the algorithms of \cite{BarWoods} do not
suffice: we actually require the list of lattice point explicitly, to
get a generating set of~$W$ as a list of monomials.  Thus the short
generating functions must be expanded.  To reduce complexity, the
short generating functions can be post-processed using the methods of
\cite{BarWoods} to yield short generating functions for the {\em
minimal}\/ generators of the ideals~$W$ in question.  Then we can
expand only these ``minimal'' short generating functions.
\end{remark}

The remaining contributions to the complexity of our local cohomology
computation come from Algorithm~\ref{a:secset}, which computes the
sets of polytopes whose disjoint unions constitute the sectors, and
Algorithm~\ref{a:secposet}.  The latter is quadratic in the output of
Algorithm~\ref{a:secset}, times a factor coming from the Minkowski sum
operations and the decision procedure for whether each such sum
contains a lattice point after intersecting with~$Q$.  Therefore it
remains only to analyze Algorithm~\ref{a:secset}.

\begin{prop} \label{p:complexity}
The number of polyhedra arising in Algorithm~\ref{a:secset} is
polynomial in the Bass numbers of~$M$ and the number of facets of~$Q$.
\end{prop}
\begin{proof}
Each Bass number of~$M$ represents an indecomposable injective module
whose bounding hyperplanes subdivide $\RR^d$ into a number of regions.
Consider the subdivision of $\RR^d$ obtained by taking simultaneously
all of the hyperplanes corresponding all of the Bass numbers of~$M$.
The number of hyperplanes contributed by each Bass number is at most
the number of facets of~$Q$, so the total number of hyperplanes is at
most the number of facets of~$Q$ times the sum of the contributing
Bass numbers.  It is well known (and follows by induction on~$n$ and
the dimension~$d$) that $n$ hyperplanes subdivide $\RR^d$ into a
number of regions that is a polynomial in~$n$ of degree~$d$.
\end{proof}

This proof shows that the number of polyhedra is exponential in the
dimension.  Exponential growth as a function of dimension also occurs
in the analysis before Remark~\ref{r:BarWoods}, where we apply
\cite{McMullen77}.

\begin{remark}\rm \label{r:BV}
A large number of rational polyhedra arise in the course of computing
local cohomology modules.  When the identification of all the lattice
points in these polyhedra is necessary, the complexity of this task
should be drastically reduced by the fact that most of these polyhedra
have facets parallel to those of~$Q$ itself.  Results such as those in
\cite{BVvectPartFctn} could be helpful along these lines.%
\end{remark}


\subsection*{Acknowledgments}
We wish to thank Gennady Lyubeznik, who motivated us to make our
previously abstract methods algorithmic.  Both authors were partially
supported by the National Science Foundation.

\footnotesize
\bigskip
\def\cprime{$'$}
\providecommand{\bysame}{\leavevmode\hbox to3em{\hrulefill}\thinspace}


\end{document}